\documentclass[12pt]{article}

\usepackage{amsmath}[1996/11/01]
\usepackage{amssymb,amsthm,amsfonts,latexsym,psfig,epsfig,here,graphics}
\usepackage{amscd, makeidx}

\def\beql#1#2\eeql{\begin{equation}\label{#1}#2\end{equation}}

\newtheorem{thm}{Theorem}
\newtheorem{theorem}[thm]{Theorem}

\newtheorem{cor}[thm]{Corollary}

\newtheorem{lemma}[thm]{Lemma}

\newtheorem{remark}[thm]{Remark}

\theoremstyle{definition}

\newtheorem{defn}[thm]{Definition}

\newcommand{\mitt}{\mbox{ : }}
\newcommand{\bew}{\noindent\underline{Proof.}\ }
\newcommand{\eb}{\phantom{zzz}\hfill{$\square $}\smallskip}

\DeclareMathOperator{\rk}{rk}
\DeclareMathOperator{\GL}{GL}
\DeclareMathOperator{\mon}{mon}
\DeclareMathOperator{\cwe}{cwe}
\DeclareMathOperator{\id}{id}
\DeclareMathOperator{\Inv}{Inv}

\DeclareMathOperator{\Aut}{Aut}

\newcommand{\C}{\mathbb C}

\newcommand{\Q}{\mathbb Q}

\newcommand{\N}{\mathbb N}
\newcommand{\F}{\mathbb F}

\begin{document}

\LARGE
\begin{center}
{\bf Kneser-Hecke-operators in coding theory.}
\end{center}

\normalsize
\begin{center}
Gabriele Nebe
\\Lehrstuhl D f\"ur Mathematik, RWTH Aachen, 52056 Aachen, Germany, nebe@math.rwth-aachen.de
\end{center}

\begin{abstract}
The Kneser-Hecke-operator is a linear operator defined on the
complex vector space spanned by the equivalence classes of a
family of self-dual codes of fixed length.
It maps a linear self-dual code $C$ over a finite field
to the formal sum of the equivalence classes of
those self-dual codes that intersect $C$ in a codimension 1 subspace.
The eigenspaces of this self-adjoint linear
operator may be described in terms of 
a coding-theory analogue of the Siegel $\Phi $-operator.
\end{abstract}

\section{Introduction}
The paper translates the lattice theoretic construction
of certain Hecke-operators from \cite{VenkovNebe}
to coding theory.
It only deals with linear self-dual codes over finite fields. 

There is a beautiful analogy between most of the notions for lattices and 
codes provided by construction A
(see for instance \cite{SPLAG}, \cite{cliff2}, \cite{Runge}).
 Theta-series of lattices correspond to weight-enumerators of codes.
Whereas theta-series of unimodular lattices are modular forms for certain
Siegel modular groups, weight-enumerators of self-dual codes
are polynomials
invariant under a certain finite group, called the associated
Clifford-Weil group, and in fact the main result of 
\cite{cliff2} shows that 
these weight-enumerators generate the invariant ring.
However, this generalized Gleason theorem
 is not true in such generality for lattices
(see for instance \cite{KoSM04}). 
One important tool in the theory of modular forms is 
Siegel's  $\Phi $-operator that maps the genus-$m$ Siegel theta-series of a 
lattice to its Siegel theta-series of genus $m-1$.
An analogue of this $\Phi $-operator was 
introduced in coding theory by Runge \cite{Runge}
to generalize Gleason's theorem to higher genus weight-enumerators of
binary self-dual codes.
Also theta-series with harmonic coefficients have a counterpart in
coding theory (see \cite{Bachoc}, \cite{Bachoc2}).
One missing concept in coding theory is that of Hecke-operators, 
which are an important tool in the theory of modular forms.
Certain of these Hecke-operators may be expressed 
in terms of lattices (see \cite{VenkovNebe}, \cite{TeiderNebe}, \cite{Krieg}).
The present paper translates this concept to coding theory.
This also answers a question raised in 1977 in \cite{Broue}.

There are slight differences from the lattice case. \\
(1) For codes this method only yields 
$p$-local Hecke-operators, where $p$  is the characteristic 
of the field, 
 whereas for lattices Hecke-operators are defined for all primes.
\\
(2)
In the theory of modular forms, the Hecke-algebra is generated by
certain double cosets of the Siegel modular group.
For codes such a commutative algebra generated by 
double cosets of the Clifford-Weil groups is 
examined in \cite{preprint}.
This maps onto the algebra generated by the Kneser-Hecke operator.
\\
(3)
The main result of the paper is that in 
the coding theory case, one can say much more about the
resulting Kneser-Hecke-operator:
The possible eigenvalues are known a priori and the 
 corresponding eigenspaces are exactly the
analogues of the spaces of Siegel cusp-forms.

\section{The general setup.}\label{gen}

Let ${\cal F}$ denote the family of self-dual codes of length $N$
of a given Type.
For a precise definition of Type  the reader is referred  to \cite{cliff2} or
\cite{cliff2a}.
The present paper only deals with codes over
finite fields $\F $ that are subspaces of $\F ^N$, self-dual with respect
to some non-degenerate bilinear or Hermitian form
$$b : \F^N \times \F^N \to \F, b(x,y) := \sum _{i=1}^N x_i \overline{y_i}$$ where
$\overline{\phantom{s}} :\F \to \F $ is either the identity (in the bilinear case)
or a non-trivial automorphism of order 2 (in the Hermitian case).
The dimension of these codes is then $n:=\frac{N}{2}$.
The integers $n$ and  $N$ will be fixed throughout the paper.

There are several possible notions of equivalence for these codes.
This paper always uses permutation equivalence,
which means that two codes are equivalent if and only if there
is a permutation $\pi \in S_N$ of the coordinates mapping one code
onto the other.
Also the automorphism group 
$$\Aut(C) = \{ \pi \in S_N \mitt \pi (C) = C \} $$
of a code is just the subgroup of the symmetric
group $S_N$ that preserves the code.
The results may be easily generalized to coarser notions of equivalence
(for instance allowing the Galois group to act)
provided that one uses the same notion for automorphisms 
of codes and that one deals with the appropriate symmetrized
 weight-enumerators  such that equivalent  codes have the same 
genus-$m$  weight-enumerator for all $m$.

Let ${\cal V}$ be the $\C $-vector space on the
set  of all equivalence classes $[C]$ with $C\in {\cal F}$.
So the set 
$${\cal B}:=\{ [C] \mitt  C \in {\cal F} \} $$
is a $\C $-basis for ${\cal V}$.

\begin{remark}
${\cal V}$ 
has 
a Hermitian positive definite scalar product defined by
$$([C],[D] ) := |\Aut(C) | \delta _{[C],[D] } $$.
\end{remark}

%

\subsection{The filtration of ${\cal V}$.}

The {\em genus-$m$ complete weight-enumerator} of a code $C$
is a homogeneous polynomial in
${\cal P}_m := \C [x_a \mitt a\in \F ^m ] $
of degree $N$.
For
 an $m$-tuple $\underline{c}:=(c^{(1)},\ldots , c^{(m)} ) \in (\F ^N)^m $ let
$$\mon (\underline{c}) = \prod _{v\in \F ^m} x_v^{a_v(\underline{c})} \in {\cal P}_m$$
where for $v = (v_1,\ldots , v_m ) \in \F ^m $ 
$$a_v(\underline{c}) = |\{ i\in \{ 1,\ldots , N \} \mitt c^{(j)}_i = v_j \mbox{ for all } 1\leq j \leq m \} |$$
is the number of columns of the $m\times N$-matrix
defined by $\underline{c}$ that are equal to $v$.
Then  $$\cwe_m (C) := \sum _{\underline{c} \in C^m} \mon (\underline{c}) 
\in {\cal P}_m $$
and $\cwe_0(C) := 1 $.
 Note that $\cwe _m(C)$ only depends on the equivalence class of $C$
and hence $\cwe_m$ may be extended to a linear map
$$\cwe _m : {\cal V} \to  {\cal P}_m, \ 
\cwe_m (\sum v_C [C] ) :=  \sum v_C \cwe _m(C)  .$$
For $m\in \N _0$ let $${\cal V}_m := \ker (\cwe _m ) \leq {\cal V}. $$

To define an analogue of 
 the Siegel $\Phi $-operator one has to choose an embedding 
$\epsilon : \F ^{m-1} \to \F ^m , 
(a_1,\ldots , a_{m-1}) \mapsto (a_1,\ldots, a_{m-1} , 0 ) .$ 
Then  there are 
 for all $m\in \N $
 ring homomorphisms 
$$\Phi  : {\cal P}_m \to {\cal P}_{m-1}
,\  x_a \mapsto \left\{ \begin{array}{ll} x_{\epsilon ^{-1}(a) } & \mbox{ if } 
a \in \epsilon (\F^{m-1}) \\ 0 & \mbox{ else. } \end{array} \right. $$
Note that $\Phi $ respects the homogeneous components of the polynomial
rings, i.e. if $p$ is a homogeneous polynomial of degree $N$ in 
${\cal P}_m $, then 
$\Phi (p) \in {\cal P}_{m-1}$ is either 0 or 
homogeneous of the same degree $N$.
Also $\Phi (\cwe_m(v) ) = \cwe _{m-1} (v) $ for all $v\in {\cal V}$.
Since the complete weight-enumerators of genus $n$ of
the basis ${\cal B}$  are
linearly independent one gets a filtration
$$ {\cal V} := {\cal V}_{-1} \geq {\cal V}_0 \geq \ldots \geq {\cal V}_{n} = \{ 0 \} $$
with $${\cal V}_0 := \{
\sum  _{ [C]\in {\cal B}} v_{C} [C]   \mitt \sum  _{ [C]\in {\cal B}} v_{C} = 0 \} $$
of codimension $1$ in ${\cal V}$.
The dual filtration is 
obtained by letting ${\cal W}_i:= {\cal V}_i^{\perp }$. Then
$$ {\cal V} = {\cal W}_{n} \geq {\cal W}_{n-1} \geq \ldots \geq {\cal W}_{0} \geq {\cal W}_{-1} = \{0 \} .$$
The space $W_0$ is one-dimensional generated by 
$${\sigma }_{N} := \sum _{[C]\in {\cal B} } |\Aut(C) |^{-1} [C] .$$

Using the Hermitian scalar product one obtains the 
orthogonal decomposition of ${\cal V}$ associated to this 
filtration by putting
$$ {\cal Y}_m :=  {\cal W}_m \cap {\cal V}_{m-1} = \{
w\in {\cal W}_m \mitt (w,x) = 0 \mbox{ for all } x\in {\cal W}_{m-1 } \} .$$
Then
\beql{decomp} 
{\cal V} = \perp _{m=0}^n {\cal Y} _m \eeql
 with
${\cal Y} _0 = {\cal W}_0 = \langle \sigma _N \rangle $.
Moreover the mapping $\cwe_m$ yields an isomorphism between
${\cal Y}_m$ and the kernel of the $\Phi $-operator 
on $\cwe_m({\cal V}) $.
One may think of the space ${\cal Y}_m$
(or the isomorphic space $\cwe _m({\cal Y}_m)$) as the analogue
of the space of Siegel cusp-forms of genus $m$.
In Section \ref{main} it is shown
that the decomposition \eqref{decomp} is in fact the eigenspace
decomposition of ${\cal V}$ under the Kneser-Hecke-operator $T$ defined below.


\subsection{Kneser-Hecke-operators.}

\begin{defn}
For $0\leq k \leq n $ 
two codes $C, D \in {\cal F}$  are called 
{\em $k$-neighbors}, written $C\sim _k D$,
 if $\dim(C\cap D) = \dim (C) - k$.
\\
Define a linear operator $T_k$ on ${\cal V}$  by
$$T_k([C]):= \sum _{D\sim_k C} [D ] $$
where the sum is over all $k$-neighbors $D \in {\cal F}$ of the code $C$.
The operator $T_k$ is called the {\em $k$-th Kneser-Hecke-operator} for ${\cal F}$.
\\
Let $T: = T_1$ be the 
{\em Kneser-Hecke-operator} and call $1$-neighbors simply {\em neighbors}.
\end{defn}

\begin{theorem}
For $0\leq k \leq n$
the operator 
$T_k$ is a self-adjoint linear operator on the vector space ${\cal V}$.
\end{theorem}

\bew
By definition $T_k$ is linear. 
For basis vectors $[C],[D] \in {\cal B}$ one has
$$\begin{array}{l}
\frac{N!}{|\Aut(D) |} |\{ C' \in {\cal F} \mitt  C'\sim _k D \mbox{ and } C'\cong C \} | \\
=
\sum _{\tilde{D} \cong D}  |\{ C' \in {\cal F} \mitt  C'\sim _k \tilde{D} \mbox{ and } C'\cong C \} | \\
=
\sum _{\tilde{C} \cong C}  |\{ D' \in {\cal F} \mitt D' \sim _k \tilde{C} \mbox{ and } D'\cong D \} | \\
=
\frac{N!}{|\Aut(C) |} |\{ D' \in {\cal F} \mitt D'\sim _k C \mbox{ and } D'\cong D \} | .\end{array} $$
The middle equality follows since the neighboring relation is symmetric 
and invariant under equivalences.
Therefore
$$\begin{array}{l}
(T_k([C]),[D] ) = 
|\Aut(D) | |\{ D' \in {\cal F} \mitt  D'\sim_k C  \mbox{ and } D'\cong D \} |  \\
=
|\Aut(C) | |\{ C' \in {\cal F} \mitt  C'\sim _k D \mbox{ and } C'\cong C \} | = ([C] ,T_k([D]) )
.\end{array} $$
Hence $T_k$ is self-adjoint.
\eb

Experiments suggest that the operators $T_k$  are polynomials
in $T=T_1$.

\subsection{The main theorem.}\label{main}

The eigenvalue of $T$ on the space ${\cal Y}_m$ 
depends on the geometry of the underlying space $(\F^N, b)$.
To prove the main theorem some more notation is needed:
Denote by
$${\cal M}_m := \{ \prod _{a\in \F ^m} x_a ^{e_a} \mitt
\sum _{a\in \F ^m} e_a = N \} \subseteq {\cal P}_m$$
the monomials in ${\cal P}_m$ of degree $N$.
\\
For a monomial $X = \prod _{a\in \F ^m} x_a ^{e_a } \in {\cal M}_m  $
define the {\em rank}
$$\rk(X) := \dim \langle a \mitt e_a > 0 \rangle $$
and let
$${\cal M}^*_m := \{ X \in {\cal M}_m \mitt \rk (X) = m \} .$$
For $X\in  {\cal M}_m $ and subset $C\subset \F ^N$ define
$$a_{X} (C) := |  \{ \underline{c} := (c^{(1)},\ldots , c^{(m)} ) \in C^m  \mitt
\mon (\underline{c} ) = X \} | . $$

\begin{remark}
(i) $a_X(C)$ only depends on the equivalence class of the code $C \leq \F^N$.
\\
(ii) $\cwe_m(C) := \sum _{X\in {\cal M}_m} a_X(C) X $.
\\
(iii)  For $X\in {\cal M}_m$ extend $a_X$ to  a linear mapping
$$a_X : {\cal V} \to \C ,\ \sum _{[C] \in {\cal B}} v_C [C] \mapsto
 \sum _{[C] \in {\cal B}} v_C a_X(C) .$$
Then 
${\cal V}_m =  \{ v\in {\cal V} \mitt a_X(v) = 0  \mbox{ for all }
X \in {\cal M}_m \} $.
\\
(iv)
${\cal V}_m = \{ v\in {\cal V}_{m-1} \mitt
a_X(v) = 0 \mbox{ for all } X \in {\cal M}^{*}_m \} .$
\end{remark}

 In this language explicit generators for the spaces
 ${\cal W}_m$ are  obtained by generalizing the construction
of $\sigma _N =b_1$.

\begin{remark}
For $X\in {\cal M}_m$ let
$$ b_X := \sum _{[C] \in {\cal B} } \frac{a_X(C)}{|\Aut(C)|}
 [C] \in {\cal V} .$$
Then for  any $v\in {\cal V}$ the scalar product
$$(b_X , v) = a_X(v) .$$
The vectors $b_X$, $X\in {\cal M}_m $ span the space ${\cal W}_m$.
\end{remark}

\bew
Let ${\cal U}_m := \langle b_X \mitt X\in {\cal M}_m \rangle $.
Then $${\cal U}_m ^{\perp } = \{
v\in {\cal V} \mitt (b_X,v) = a_X(v) = 0 \mbox{ for all } X\in {\cal M} _m \} =
{\cal V}_m = {\cal W}_m^{\perp } $$
and therefore ${\cal U}_m = {\cal W}_m$.
\eb

\begin{remark}{\label{basisunab}}
Let $\underline{c} := (c^{(1)}, \ldots , c^{(m)} ) \in C^m $,
$X:= \mon(\underline{c} )$,
and $U:= \langle c^{(1)}, \ldots , c^{(m)} \rangle \leq C $.
Then $\dim (U) = \rk (X) $.
If $\rk (X) = m $ 
 and
$\underline{b} := (b^{(1)}, \ldots , b^{(m)} ) $ is another basis of $U$,
then for all $v\in {\cal V}$ $$a_X(v) = a_{\mon(\underline{b})} (v) .$$
\end{remark}

For the proof of the main theorem 
 choose a suitable subset 
${\cal M}^0_m \subset {\cal M}^*_m $ such that 
\beql{cond}
 {\cal V}_m = \{ v\in {\cal V}_{m-1} \mitt 
a_X(v) = 0 \mbox{ for all } X \in {\cal M}^{0}_m \} . 
\eeql
Clearly the full set ${\cal M}^0_m = {\cal M}^* _m$ satisfies
condition \eqref{cond},  but there may be smaller sets.

\begin{lemma}\label{m1}
Assume that all codes in ${\cal F}$ contain the 
 all-ones vector ${\bf 1}:= (1,\ldots ,1) $.
Then
$${\cal M}^0_m = {\cal M}^1_m := 
\{ \mon (\underline{c}) \mitt \dim \langle {\bf 1}, 
c^{(1)}, \ldots , c^{(m)} \rangle  = m+1 \}  $$
satisfies condition \eqref{cond}.
\end{lemma}
\bew
For $C\in {\cal F}$
 let $\underline{c} = (c^{(1)},\ldots , c^{(m)} ) \in C^m $
be such that $X:=\mon (\underline{c} ) \in {\cal M}_m^* \setminus {\cal M}_m^1$.
Then ${\bf 1} \in U := \langle c^{(1)},\ldots , c^{(m)} \rangle $ and
$\dim (U) = m$.
By Remark \ref{basisunab}
 $a_X([C] )$ is independent of the choice of the basis $\underline{c}$ of
$U$  one may assume w.l.o.g. that $c^{(m)} = {\bf 1}$.
Let
  $\underline{c'} := (c^{(1)},\ldots , c^{(m-1)} ) $ and
$Y:=\mon (\underline{c'})  \in {\cal M}_{m-1} ^*$.
Then $$\underline{b'} := (b^{(1)},\ldots , b^{(m-1)})  \mapsto
\underline{b} := (b^{(1)},\ldots , b^{(m-1)}, {\bf 1} ) $$
establishes a bijection between
$$\{ \underline{b'} \in C ^{m-1}  \mitt \mon (\underline{b'} ) = Y \}
\mbox{ and }
\{ \underline{b} \in C ^{m}  \mitt \mon (\underline{b} ) = X \}  $$
showing that  $a_X([C]) = a_Y([C]) $  for all $[C]\in {\cal B}$.
Therefore $a_X(v) = 0 $ for all $v\in {\cal V}_{m-1} \subseteq
\ker (a_Y )$.
Hence
$$
\begin{array}{ll}
{\cal V}_m  & = \{ v\in {\cal V}_{m-1} \mitt a_X(v) = 0 \mbox{ for all } X\in {\cal M}_m^* \} \\
 & = \{ v\in {\cal V}_{m-1} \mitt a_X(v) = 0 \mbox{ for all } X\in {\cal M}_m^1 \}
\\ \end{array}
$$
which shows condition \eqref{cond}.
\eb

{\bf Condition $\star $.}
{\rm In addition to condition 
\eqref{cond},  assume that 
 for all codes $C\in {\cal F}$ and 
all  $\underline{c} := (c^{(1)}, \ldots , c^{(m)}) \in C^m$ such that
$\mon (\underline{c} ) \in {\cal M}^0_m$ the sum 
$$ \alpha _m := \sum _{E\in {\cal E}_C(\underline{c})}  \alpha _E $$  
does not depend on $\underline{c} $ and $C$.
Here $${\cal E}_C(\underline{c}) := \{ E\leq C \mitt \dim (E) = n-1, \underline{c} \in E^m \}  $$
and  $\alpha _E = \alpha _E (C)$ is the number of codes $D\in {\cal F}$ with 
$D\cap C = E$.}
\\

Furthermore let 
$$\beta _m := \frac{|\F|^m-1}{|\F|-1}  $$
the number of $(m-1)$-dimensional subspaces of $\F ^m$
and, if Condition $\star $ is satisfied, 
$$\nu _m := \alpha _m  - \beta_m .$$

\begin{theorem}{\label{tmain}}
Assume that Condition $\star  $ is satisfied. 
Then the
 space ${\cal Y}_m $ is exactly the $\nu _m$-eigenspace of $T$ in ${\cal V}$.
\end{theorem}

\bew
It is enough to show that $T$ acts as 
$\nu _m \id $  on ${\cal V}_{m-1}/{\cal V}_m$ which just means that 
for $$v := \sum _{[C]\in {\cal B} } v_C [C] \in {\cal V}_{m-1} = \ker (\cwe _{m-1} ) $$
the difference
$$T(v) - \nu _m v \in {\cal V}_{m} = \ker (\cwe _m ) .$$
For this it is enough to show that 
$$a_X(T(v)) = \nu _m a_X (v) \mbox{ for all } X\in {\cal M}_m^0 .$$
Now $$T(v) = \sum _{[C] \in {\cal B} } v_C  T([C]) = \\
 \sum _{[C] \in {\cal B} } v_C   
\sum _{\small \begin{array}{c} E\leq C \\ \dim (E) = n-1\end{array}} \sum _{\small
\begin{array}{c} D\in {\cal F} \\ E= D\cap C \end{array} } [D] $$
therefore we have to calculate for $X\in {\cal M}_m^{0} $ and 
a fixed $C\in {\cal F} $
\beql{sum}
\sum _{\small \begin{array}{c} E\leq C \\ \dim (E) = n-1 \end{array}} \sum _{\small \begin{array}{c} D\in {\cal F} \\ E= D\cap C  \end{array}} a_X ( [D]).
\eeql
\\
Let
 $\underline{c}:=(c^{(1)}, \ldots , c^{(m)})\in D^m$ for some neighbor
 $D$  of $C$ 
such that $\mon (\underline{c}) = X \in {\cal M}^0_m$.
Put $W:= \langle c^{(1)}, \ldots , c^{(m)}\rangle $ and distinguish two cases:
\\
(a) $W\leq C$:
Then $\underline{c} \in D^m$, if and only if $\underline{c} \in (D\cap C)^m$ and
by Condition $\star $
 this yields a contribution $\alpha _m a_X(C) $ to the
sum \eqref{sum}.
\\
(b)  $W\not\leq C$:
Then 
 $U:=W \cap C $ has dimension $m-1$, and 
$E = W ^{\perp }  \cap C $ and 
$D = \langle C , W \rangle$ 
are uniquely determined by $\underline{c} $. 
Let $\underline{b}:=( b^{(1)}, \ldots , b^{(m-1)} ) $ be a basis of $U$ and let 
$Y:= \mon (\underline{b}) $. Then by Remark \ref{basisunab}
the value $a_Y (C) $ is independent of the
choice of this basis.
Note that $W $ has exactly $\beta _m$ 
such submodules $U$.
Here the contribution to the sum \eqref{sum} 
is $$\sum _{Y \leq X} (a_Y (C) - a_X(C) ) = 
(\sum _{Y \leq X} a_Y (C)) - \beta _m a_X(C) ~. $$
By induction on $m$ this argument shows that the subspaces 
${\cal V}_{m}$ are $T$-invariant. 
Furthermore for $v\in \ker(\cwe_{m-1} )$ the sum
$$a_Y(v) = \sum _{[C]\in {\cal B}} v_C a_Y([C]) = 0 \mbox{ for all } Y\in {\cal M}_{m-1} .$$
Hence in total 
$$ a_X(T(v)) = 
\nu _m a_X(v) $$
 for all $X\in {\cal M}_m^0 $ and $v\in {\cal V}_{m-1}$.
\eb

%

%
%

\section{Classical Types.}

In this section it is shown that all classical Types of self-dual 
codes over finite fields  $\F = \F_q$
satisfy Condition $\star $ of Section \ref{main} and the eigenvalues of 
the operator $T$ are determined.
The Types are denoted by the names used 
in \cite{chapter} and \cite{cliff2}:

\begin{itemize}
\item[$q^E$]:
Euclidean self-dual $\F_q$-linear codes in odd characteristic.
So ${\cal F} = \{ C = C^{\perp } \leq \F_q^N \} $ where the dual code
$C^{\perp} = \{ v\in \F_q^N \mitt b(v,c) := \sum _{i=1}^N v_ic_i = 0
\mbox{ for all } (c_1,\ldots , c_N )\in C \} $.
 \item[$q^E_{\bf 1}$]:
 Same as $q^E$ but we additionally impose the condition that
 the all-ones vector ${\bf 1} = (1,\ldots , 1)$ be in all codes 
in ${\cal F}$.
 \item[$q^E_{\rm I}$]:
 Same as $q^E$  but now  $q$ is even.
 \item[$q^E_{\rm II}$]:
Same as $q^E_{\rm I}$ but additionally assuming that
 the codes in ${\cal F}$ 
are generalized doubly even as defined in \cite{quebbemann}, \cite{NQRS}.
 \item[$q^H$]:
 Hermitian self-dual $\F_q$-linear codes.
 Here $q=r^2$ is a square and $\overline{\phantom{s}} : \F_q \to \F_q, x \mapsto x^r $ denotes the non-trivial Galois automorphism of $\F _q/\F_r $.
 Then the dual code 
$ C^{\perp} := \{ v\in \F_q^N \mitt \sum _{i=1}^N c_i\overline{v_i} = 0 \mbox{
for all } c\in C \} \}
$
\item[$q^H_{\bf 1}$]:
 Same as $q^H$, but additionally assuming that $\bf 1$
 be in the codes in ${\cal F}$.
\end{itemize}

To show that these Types satisfy condition $\star $ from 
Section \ref{main} we need to precise the set 
${\cal M}_m^{0}$ and calculate the number 
$\alpha _m$ as defined there. 
Then the eigenvalue of $T$ on ${\cal Y}_m$ is 
$\nu _m = \alpha _m - (q^m-1)/(q-1) $ according to Theorem \ref{tmain}.

\begin{theorem}
The codes of the six Types listed above satisfy condition $\star $.
The following table which lists the sets ${\cal M}_m^0$,
the corresponding value for
$\alpha _m$ and the eigenvalue $\nu _m$  (multiplied by $q-1$ to avoid 
fractions):
$$
\begin{array}{|c|c|c|c|}
\hline
\mbox{Type } & {\cal M}_m^0 & \alpha _m (q-1) & \nu _m (q-1) \\
\hline
\hline
q^E_{{\rm I}} & {\cal M}_m^1 & q^{n-m}-q & q^{n-m}-q-q^m+1 \\
\hline
q^E_{\rm II} & {\cal M}_m^1 &  q^{n-m-1}-1 & q^{n-m-1}-q^m \\
\hline
q^E & {\cal M}_m^* & q^{n-m}-1 & q^{n-m}-q^m \\
\hline
q^E_{1} & {\cal M}_m^1 & q^{n-m-1}-1 & q^{n-m-1}-q^m \\
\hline
q^H & {\cal M}_m^* & \sqrt{q}(q^{n-m}-1) & q^{n-m+1/2}-q^m-q^{1/2}+1 \\
\hline
q^H_1 & {\cal M}_m^1 & \sqrt{q}(q^{n-m-1}-1) & q^{n-m-1/2}-q^m-q^{1/2}+1 \\
\hline
\end{array}
$$
\end{theorem}

\bew
By Lemma \ref{m1} the set
${\cal M}_m^1$ satisfies the condition 
\eqref{cond} of Section \ref{main} in the cases where  the codes in
${\cal F}$ contain the all-ones vector.
\\
Let $C\in {\cal F}$ be a self-dual code of one of the six Types and 
let $E\leq C$ be a subspace of codimension 1.
First the number $\alpha _E := | \{ D\in {\cal F} \mid C\cap D = E \} |$
is determined.
The relevant codes $D$ correspond to the
one-dimensional isotropic subspaces $\neq C/E$ of 
$E^{\perp }/E $ with respect to the associated geometry.
If ${\bf 1} \in D$ for all $D\in {\cal F}$ (which is the case for 
$q^E_{\rm I}, q^E_{\rm II}, q^E_{1}, q^H_1$) and ${\bf 1}\not\in E$ then 
$C = \langle E,{\bf 1} \rangle $ is the unique code in ${\cal F}$ that
contains $E$.
So here $\alpha _E = 0$ if ${\bf 1}\not\in E$.

\underline{Case $q^E_{\rm I}$}:
Assume that we are in case $q^E_{\rm I}$ and that ${\bf 1}\in E$.
Then all elements $c=(c_1,\ldots, c_N) \in E^{\perp } $
satisfy
$$ 0 = b( {\bf 1} , c) = \sum _{i=1}^N  c_i
 = \sum _{i=1}^N c_i^2
 = (\sum _{i=1}^N c_i ) ^2
= b(c,c) ,$$
because the characteristic of $\F $ is 2.
Hence all $q+1$ one-dimensional subspaces of $E^{\perp }/ E$ are
self-dual and $\alpha _E = q$.
This proves that for 
$\underline{c} = (c^{(1)},\ldots , c^{(m)} ) \in C^m $
with $X:=\mon (\underline{c} ) \in {\cal M}_m^1$ the sum 
$ \alpha _m:=\sum _{E\in {\cal E}_C(\underline{c})}  \alpha _E $ is $q$ times 
the number of 
$(n-1)$-dimensional subspaces of $C$ that contain the $(m+1)$-dimensional
space $\langle {\bf 1},c^{(1)},\ldots ,c^{(m)} \rangle $
hence $\alpha _m = q \beta _{n-m-1}= q \frac{q^{n-m-1}-1}{q-1} $.

\underline{Case $q^E_{\rm II}$}:
If ${\bf 1}\in E$, the space $E^{\perp }/ E$ is 
a non-singular quadratic space 
of dimension 2 with a maximal isotropic
subspace $C/E$.
Hence $E^{\perp}/E $ is a hyperbolic plane and has exactly 
two maximal isotropic subspaces.
Therefore $\alpha _m =\beta_{n-m-1} = \frac{q^{n-m-1}-1}{q-1} $.

\underline{Case $q^E$}:
Here $q$ is odd and for any codimension 1 subspace $E\leq C\in {\cal F}$
the space $E^{\perp }/E$ is a hyperbolic plane with exactly 
two maximal isotropic subspaces.
Therefore $\alpha _E = 1$ for all $E$ and 
$\alpha _m = \beta _{n-m} = \frac{q^{n-m}-1}{q-1} $.

\underline{Case $q^E_1$}:
By the argumentation above we only need to consider the subspaces 
$E$ that contain ${\bf 1}$.
Again $E^{\perp}/E$ is a hyperbolic plane and
$\alpha _m = \beta _{n-m-1} = \frac{q^{n-m-1}-1}{q-1} $.

\underline{Case $q^H$}:
Let $E \leq C \in {\cal F}$ be a self-orthogonal
subspace of dimension $n-1$.
Then $E^{\perp }/ E $ is a 2-dimensional non-degenerate
Hermitian space over $\F _q $ hence isometric to
$\F _q ^2 $ with the Hermitian form $((x_1,x_2),(y_1,y_2)) := x_1\overline{y_1}
+ x_2 \overline{y_2} $.
It follows easily that this space has exactly $\sqrt{q}+1$ one-dimensional
isotropic subspaces.
Since $C/E $ is one of them $\alpha _E = r= \sqrt{q} $ and 
$\alpha _m = \sqrt{q} \beta _{n-m} = \sqrt{q} \frac{q^{n-m}-1}{q-1} $.

\underline{Case $q^H_1$}:
Similar as for $q^H$ but we only need to consider the subspaces
$E$ that contain ${\bf 1}$.
Therefore 
$\alpha _m = \sqrt{q} \beta _{n-m-1} = \sqrt{q} \frac{q^{n-m-1}-1}{q-1} $.
\\
The eigenvalue $\nu _m$ of $T$ now results from the general formula 
in Theorem \ref{tmain}.  \eb
\subsection{Explicit numerical results.}

The neighboring method provides a quite efficient way to enumerate
all equivalence classes of codes of a given Type.
During this procedure, the Kneser-Hecke-operator $T$ is calculated
without difficulty.
It is then easy to obtain the eigenvalues of $T$ and
the (dimensions of the) eigenspaces.
With Theorem \ref{tmain} this gives the dimension of all 
spaces $\cwe _m({\cal V})$ for all $m\in \N _{0 }$,
even if it might be quite difficult to obtain 
(enough terms of) the genus-$m$ complete weight-enumerators 
of the codes in ${\cal F}$ to calculate this dimension 
directly.
The calculations are performed with MAGMA \cite{magma} using
a direct analogue of the Kneser-neighboring procedure 
described in \cite{Kneser} for lattices.
Recall that  $n :=\frac{N}{2}$ denotes the dimension of 
the codes in ${\cal F}$.
Starting with some code in $C\in {\cal F}$ (usually constructed 
as an orthogonal sum) we enumerate the orbits of the
automorphism group on the $(n-1)$-dimensional subspaces $E\leq C$
(resp. those $E$ that contain ${\bf 1}$) and calculate the 
neighbors of $C$ as preimages of the isotropic one-dimensional
subspaces of $E^{\perp }/ E$.

For the binary codes (where calculations could be performed 
without problems up to  length $N=32$) we have the following 
explicit results.

\begin{center}
{\bf Table 1: The dimension of the space ${\cal Y}_m$ for Type $2_{\rm I}$.}
\end{center}
$$
\begin{array}{|c|c|c|c|c|c|c|c|c|c|c|c|c|}
\hline
N, m & 0 & 1 & 2 & 3 & 4 & 5 & 6 & 7 & 8 & 9 & 10 & 11   \\
\hline
2& 1 &  &  &  &  &  &  &  &  &  &  &     \\
4& 1 &  &  &  &  &  &  &  &  &  &  &     \\
6& 1 &  &  &  &  &  &  &  &  &  &  &     \\
8& 1 & 1&  &  &  &  &  &  &  &  &  &     \\
10& 1 &1 &  &  &  &  &  &  &  &  &  &     \\
12& 1 &1 & 1&  &  &  &  &  &  &  &  &     \\
14& 1 & 1 & 1 & 1 &  &  &  &  &   &  &    &   \\
16& 1 & 2& 1& 2& 1&  &  &  &  &  &  &     \\
18& 1 & 2& 2&2 & 2&  &  &  &  &  &  &     \\
20& 1 & 2& 3& 4& 4& 2&  &  &  &  &  &     \\
22& 1 & 2& 3& 6& 7& 4& 2&  &  &  &  &     \\
24& 1 & 3& 5& 9&15 & 13 & 7 &2  &  &  &  &     \\
26& 1 & 3 & 6 & 12 & 23 & 29 & 20 & 8 & 1 &  &    &   \\
28& 1 & 3& 7&18&40&67&75&39&10&1 &  &     \\
30& 1 & 3 & 8 & 23 &65  & 142 & 228 & 189 & 61 & 10 & 1 &     \\
32& 1 & 4 & 10 & 33 & 111 & 341 & 825 & 1176 & 651 & 127 & 15 & 1    \\
\hline
\end{array}
$$

\begin{center}
{\bf Table 2: The dimension of the space ${\cal Y}_m$ for Type $2_{\rm II}$.}
$$
\begin{array}{|c|c|c|c|c|c|c|c|c|c|c|c|}
\hline
N, m & 0 & 1 & 2 & 3 & 4 & 5 & 6 & 7 & 8 & 9 & 10   \\
\hline
8& 1 & &  &  &  &  &  &  &  &  &     \\
\hline
16& 1 & 0& 0& 1& &  &  &  &  &  &     \\
\hline
24& 1 & 1& 1& 2& 2& 1& 1&  &  &  &     \\
\hline
32& 1 & 1& 2&5 & 10& 15& 21 & 18& 8&3 &1    \\
\hline
\end{array}
$$
\end{center}

{\bf Application to Molien-series.}

By \cite{Runge} (see also \cite{cliff1}) there is a finite
matrix group ${\cal C}_m \leq \GL _{2^m}(\Q[\sqrt{2}] )$ called the
{\em real Clifford-group of genus $m$}, such that the invariant
ring of ${\cal C}_m$ is the image of $\cwe _m $,
$$\Inv ({\cal C}_m) = \bigoplus _{N=0}^{\infty }
\langle \cwe _m (C) \mitt C=C^{\perp } \leq \F_2^N \rangle ~.$$

\begin{cor}
For $m\ge 1 $
the Molien series of ${\cal C}_m$ is
$$1+ t^2 + t^4 + t^6 +2t^8 +2t^{10 } + \sum _{N=12}^{\infty } a_N (m) t^N
$$
where for $N\leq 32$ the coefficients
$a_N(m) = \dim \langle \cwe _m(C) \mitt C=C^{\perp} \leq \F _2^N \rangle$
 are given in the following
table.
\end{cor}

$$ \begin{array}{|c|c|c|c|c|c|c|c|c|c|c|c|}
\hline
N & 12 & 14 & 16 & 18 & 20 & 22 & 24 & 26 &28 & 30 & 32  \\
\hline
m = 1 & 2 & 2 & 3 & 3 & 3 & 3 & 4 & 4 & 4 & 4 & 5 \\
\hline
m = 2 & 3 & 3 & 4 & 5 & 6 & 6 & 9 & 10 & 11 & 12 & 15 \\
\hline
m = 3 & 3 & 4 & 6& 7 & 10& 12 & 18 & 22 & 29 & 35 & 48 \\
\hline
m=4 & 3 & 4 & 7 & 9 & 14 & 19 & 33 & 45 & 69 & 100 & 159 \\
\hline
m=5 & 3 & 4 & 7 & 9 & 16 & 23 & 46 & 74 & 136 & 242 & 500 \\
\hline
m=6 & 3 & 4 & 7 & 9 & 16 & 25 & 53 & 94 & 211 & 470 & 1325 \\
\hline
m=7 & 3 & 4 & 7 & 9 & 16 & 25 & 55 & 102 & 250 & 659 & 2501 \\
\hline
m= 8  & 3 & 4 & 7 & 9 & 16 & 25 & 55 & 103 & 260 & 720 & 3152 \\
\hline
m=9  & 3 & 4 & 7 & 9 & 16 & 25 & 55 & 103 & 261 & 730 & 3279  \\
\hline
m=10  & 3 & 4 & 7 & 9 & 16 & 25 & 55 & 103 & 261 & 731 & 3294  \\
\hline
m \geq 11  & 3 & 4 & 7 & 9 & 16 & 25 & 55 & 103 & 261 & 731 & 3295 \\
\hline
\end{array}
$$

Similarly the genus-$m$ complete weight-enumerators of the
doubly-even self-dual binary codes span the invariant ring of the
{\em complex Clifford-group} ${\cal X}_m \leq \GL _{2^m} (\Q [\zeta _8])$
(see \cite{Runge}, \cite{cliff1}) and Table 2 above gives the
first terms of the Molien series of those groups.
The full Molien series of ${\cal C}_m$ and ${\cal X}_m$  are known
for $m\leq 4$ (see \cite{Our97},
sequences number A008621, A008718, A024186, A110160,
 A008620, A028288, A039946, A051354 in \cite{OEIS}).

\end{document}